\documentclass{article}
\usepackage{amsmath,amsthm,amssymb}
\usepackage[all]{xy}

\newtheorem{theor}{Theorem}
\newtheorem{propo}{Proposition}
\newtheorem{lem}{Lemma}
\newtheorem{coro}{Corollary}
\theoremstyle{definition}
\newtheorem{defi}{Definition}

\newcommand{\Z}{\mathbb Z}
\newcommand{\V}{{\mathcal V}}
\newcommand{\wt}{\widetilde}
\newcommand{\cor}{\mathrm{corank\,} }
\newcommand{\rank}{\mathrm{rank\,} }
\newcommand{\sgn}{\mathrm{sgn}}
\newcommand{\fr}{\mathrm{fr}}
\newcommand{\Kh}{\mathrm{\overline{Kh}}}

\newcommand{\tar}[2]{\begin{array}{c} #1 \\ \mbox{#2}\end{array}}
\newcommand{\bega}{\left(\begin{array}}
\newcommand{\ena}{\end{array}\right)}

\title{Khovanov homology of graph-links}
\author{Igor Nikonov}
\date{}

\begin{document}
\maketitle

\begin{abstract}
We define Khovanov homology over $\Z_2$ for graph-links.
\end{abstract}

\section{Introduction}

Graph-links~\cite{IM,IM1} are combinatorial analog of classical and
virtual links. "Diagrams" of graph-links are simple undirected
labeled graphs and graph-links themself are equivalence classes of
the graphs modulo formal Reidemeister moves. In knot theory
graph-links appear as intersection graphs of rotating chord diagrams
of links. It is known that intersection graph determines chord
diagram uniquely up to mutations~\cite{CL}. Thus any link invariant
that does not distinguish mutant links is a candidate for an
invariant of graph-links. For example, Alexander polynomial, Jones
polynomial and HOMFLY are of such type. The problem is the
equivalence relation of the graph-links can be too strong so there
can exist two chord diagrams which have different values of the
invariant but whose intersection graphs coincide as graph-links. So
the problem of finding a natural construction that would extend 
the invariant from intersection graphs to graph-links is nontrivial.
For Jones polynomial this problem was solved by D.~Ilyutko and
V.~Manturov~\cite{IM}. This article provides a solution for (odd)
Khovanov homology with coefficients in $\Z_2$.

A theory analogous to the theory of graph-links was developed by
L.~Traldi and L.~Zulli~\cite{TZ}. Their approach (called looped
interlacement graphs theory)  is based on Gauss diagrams of knots
instead of rotating chord diagrams. D.~Ilyutko showed these two
theories are closely related~\cite{I}. So the construction of the
Khovanov homology described below will provide as well (after the
necessary reformulation) a definition of the Khovanov homology for
the looped interlacement graphs.

Odd Khovanov homology of links was defined by P.~Ozsv\'ath,
J.~Rasmussen and Z.~Szab\'o~\cite{ORS}. J.~Bloom proved that odd
Khovanov homology is mutant invariant~\cite{Blo} and we can expect
it to be extended to an invariant of graph-links. In fact the
construction of J.~Bloom over $\Z_2$ can be directly tranfered to the
graph-link theory. One can also consider this invariant as the
ordinary Khovanov homology of graph-links because for links both
homology theories coincide modulo 2. On the other hand,
integer-valued Khovanov homology is not mutant invariant~\cite{Weh}
so we don't anticipate any adaption of it for graph-links.

\section{Graph-links}

Let G be an undirected graph without loops and multiple edges and
$\V=\V(G)$ be the set of its vertices. We assume $G$ be a {\em
labeled} graph, i.e. every vertex $v$ in $G$ is endowed with a pair
$(a, \alpha)$ where $a\in\{0,1\}$ is the {\em framing} $\fr(v)$ of $v$
and $\alpha \in\{-1, 1\}$ is the {\em sign} $\sgn(v)$ of the vertex
$v$.

Let us fix an enumeration of vertices of $G$. We define the {\em
adjacency matrix} $A(G)=(a_{ij})_{i,j=1,\dots,n}$ over $\Z_2$ as
follows: $a_{ij}=1$ if and only if the vertices $v_i$ and $v_j$ are
adjacent, and $a_{ij}=0$ otherwise. Besides we set
$a_{ii}=\fr(v_i)$.

Let $v\in \V$. The set of all vertices in $\V$ adjacent to $v$ is
called the {\em neighbourhood} of the vertex $v$ and denoted
$N(v)$.

Let us define the Reidemeister moves of labeled graphs.

$\Omega_1$. The first Reidemeister move is an addition/removal of an
isolated vertex labeled $(0,\pm 1)$.

$\Omega_2$. The second Reidemeister move is an addition/removal of
two adjacent (resp. nonadjacent) vertices having the labels $(1,-1)$
and $(1,1)$ (resp. $(0,-1)$ and $(0,1)$) and the same
neighbourhoods.

$\Omega_3$. The third Reidemeister move is defined as follows. Let
$u, v, w$ be three vertices of $G$ having the labels $(0,-1)$ and
$u$ be adjacent only to $v$ and $w$. Then we only change the
adjacency of $u$ with the vertices $t\in N(v)\setminus N(w)\cup
N(w)\setminus N(v)$ (for other pairs of vertices we don't change
their adjacency). In addition, we change the signs of $v$ and $w$ to
'+1'. The inverse operation is also called the third Reidemeister
move.

$\Omega_4$. The fourth Reidemeister move is defined as follows. We
take two adjacent vertices $u$ labeled $(0, \alpha)$ and $v$ labeled
$(0,\beta)$. Then we change the label of $u$ to $(0, -\beta)$ and
the label of $v$ to $(0, -\alpha)$. Besides, we change the adjacency
for each pair $(t,w)$ of vertices where $t\in N(u)$ and $w\in
N(v)\setminus N(u)$ or $t\in N(v)$ and $w\in N(u)\setminus N(v)$.

$\Omega'_4$. In this fourth move we take a vertex $v$ with the label
$(1,\alpha)$. We change the adjacency for each pair $(t,w)$ of
vertices where $t, w\in N(u)$. Besides, we change the sign of $v$
and the framing for each $w\in N(u)$.

\begin{defi}
A {\em graph-link} is the equivalence class of a simple labeled
graph modulo $\Omega_1-\Omega'_4$ moves.
\end{defi}

\section{Khovanov homology of graph links}

Let $G$ be a simple labeled graph with $n$ vertices and $A=A(G)$ be
its adjacency matrix. We now formulate mod $2$ version of
Bloom's construction of odd Khovanov homology. The reason why we
use $\Z_2$ is the chain complex construction depends on the
coefficients of the adjacency matrix $A$ and the matrix $A$ is
defined over $\Z_2$.

Suppose $s\subset \V=\V(G)$. We shall call subsets of $\V$ as {\em
states}. Define $G(s)$ to be the complete subgraph in $G$ with the
set of vertices $s$ and denote $A(s)=A(G(s))$. Consider the vector
space
$$ V(s) = \Z_2<x_1, \dots, x_n\,|\, r^s_1,\dots, r^s_n>$$ where
relations $r^s_1,\dots, r^s_n$ are given by the formula
\begin{equation}
r^s_i=\left\{\begin{array}{cl} x_i+\sum\limits_{\{j\,|\,v_j\in s\}}
a_{ij}x_j,& \mbox{if }v_i\not\in s,\\
\sum\limits_{\{j\,|\,v_j\in s\}} a_{ij}x_j,& \mbox{if }v_i\in s
 \end{array}\right.
\end{equation}
The dimension of $V(s)$ is equal to $\cor A(s)$.

There is a natural bijection between states $s\subset\V$ and
vertices of the hypercube $\{0,1\}^n$: the state $s$ corresponds to
the vector $(\alpha_1,\dots,\alpha_n)$ where $\alpha_i=0$ if $v_i\in
s$ and $\sgn(v_i)=1$ or if $v_i\not\in s$ and $\sgn(v_i)=-1$ and
$\alpha=1$ otherwise. Every edge of the hypercube is of the type
$s\rightarrow s\oplus i$ where $s\oplus i$ denotes $s\cup\{v_i\}$ if
$v_i\not\in s$ and $s\setminus\{v_i\}$ if $v_i\in s$. We orient the
arrow so that $v_i\not\in s$ if $\sgn(v_i)=-1$ and $v_i\in s$ if
$\sgn(v_i)=1$.

We assign to every edge $s\rightarrow s\oplus i$ the map
$\partial_{s\oplus i}^s : \bigwedge^*V(s)\to\bigwedge^*V(s\oplus i)$
of exterior algebras defined by the formula
\begin{equation}
\partial_{s\oplus i}^s(u)=\left\{\begin{array}{cl} x_i\wedge u& \mbox{if }x_i=0\in V(s),\\
u& \mbox{if }x_i\ne 0\in V(s).
 \end{array}\right.
\end{equation}

Consider the chain complex
$$C(G)=\bigoplus\limits_{s\subset\V}\bigwedge\nolimits^*V(s)$$
with differential
$$\partial(u)=\sum\limits_{\{s,s'\subset\V\,|\, s\rightarrow s'\}}
\partial_{s'}^s(u).$$

\begin{propo}\label{prop_val}
Chain complex $(C(G),\partial)$ is well defined.
\end{propo}
\begin{proof}%[Proof of Proposition~\ref{prop_val}]
We need to show that the maps $\partial_{s\oplus i}^s$ are well
defined and that each 2-face in the hypercube of states is commutative.

\begin{lem}\label{lem1} Let us consider a state $s$ and an index $i$ such that $v_i\not\in
s$. We can assume without loss of generality  that $A(s\oplus i) =
\bega{cc} A& a^\top\\ a&\alpha\ena$ where $A=A(s)$. Then
\begin{enumerate}
\item $x_i=0\in V(s)$ iff $\rank A=\rank \bega{c} A\\ a\ena$;
\item $x_i=0\in V(s\oplus i)$ iff $\rank \bega{c} A\\ a\ena+1=\rank \bega{cc} A& a^\top\\
a&\alpha\ena$.
\end{enumerate}
\end{lem}

\begin{proof}
It follows from the relation $x_i =\sum\limits_{\{j\,|\,v_j\in s\}}
a_{ij}x_j$ that the equality $x_i=0\in V(s)$ means the vector $a$
depends on the rows of the adjacency matrix $A$. This leads to the
first statement of the lemma.

The equality $x_i=0\in V(s\oplus i)$ implies that the vector $(0\
1)$ depends on the rows of the matrix $A(s\oplus i)$. Then
$\rank \bega{cc} A& a^\top\\
a&\alpha\ena = \rank \bega{cc} A& a^\top\\
a&\alpha\\ 0 &1\ena$. But
$$
\rank \bega{cc} A& a^\top\\
a &\alpha\\ 0 &1\ena = \rank \bega{cc} A& 0\\
a&0\\ 0 &1\ena = \rank \bega{c} A\\
a\ena+1.
$$
\end{proof}

\begin{lem}\label{lem2}
For any symmetric matrix $A$ if $\rank \bega{cc} A& a^\top\\
a&\alpha\ena=\rank A+1$ then $\rank \bega{c} A\\
a\ena=\rank A$.
\end{lem}
\begin{proof}
Assume that $\rank \bega{c} A\\
a\ena\ne \rank A$. Then $\rank \bega{cc} A& a^\top\\
a&\alpha\ena = \rank \bega{c} A\\
a\ena$ and the vector $\bega{c} a^\top\\
\alpha\ena$ depends on the columns of the matrix $\bega{c} A\\
a\ena$. Then the vector $a^\top$ depends of the columns of $A$ and
(after transposition) the vector $a$ depends on the rows of $A$.
Hence, $\rank A = \rank\bega{c} A\\ a\ena$ that contradicts to the
initial assumption.
\end{proof}

\begin{coro} For any state $s$ and any index $i$ we have
\begin{enumerate}
\item $\dim V(s\oplus i)=\dim V(s)+1$ iff $x_i=0\in V(s)$ and $x_i\ne0\in V(s\oplus
i)$;
\item $\dim V(s\oplus i)=\dim V(s)-1$ iff $x_i\ne 0\in V(s)$ and $x_i=0\in V(s\oplus
i)$;
\item $\dim V(s\oplus i)=\dim V(s)$ iff $x_i=0\in V(s)$ and $x_i=0\in V(s\oplus
i)$.
\end{enumerate}
The case $x_i\ne 0\in V(s)$ and $x_i\ne 0\in V(s\oplus i)$ is
impossible.
\end{coro}
\begin{proof}
The statements of the corollary follow from
Lemmas~\ref{lem1},\ref{lem2} and the fact $\dim V(s)=\cor A(s)$,
$\dim V(s\oplus i)=\cor A(s\oplus i)$.
\end{proof}

We call the first two cases in the proposition {\em even} and the
third case {\em odd}. From the definition of the differentials we
have $\partial^s_{s\oplus i}=0$ in the odd case.

\noindent{\em Correctness of chain maps.}

Let us consider the map $\partial_{s\oplus i}^s : \bigwedge^*V(s)\to
\bigwedge^*V(s\oplus i)$. We must check that the relations are mapped to
relations, that is for any element $u$ and any index $j$ there exist
elements $u_k\in V(s\oplus i)$ such that
$$\partial_{s\oplus i}^s(r^s_j\wedge u) = \sum_k
r^{s\oplus i}_k\wedge u_k\ \in V(s\oplus i).$$

For any $j$ we have $r_j^s = r_j^{s\oplus i} +\alpha_j x_i$ for some
$\alpha_j\in\Z_2$. If $x_i=0\in V(s\oplus i)$ then
$$\partial_{s\oplus
i}^s(r^s_j\wedge u)=r^s_j\wedge u' = r^{s\oplus i}_j\wedge u'+\alpha_j
x_i\wedge u'=r^{s\oplus i}_j\wedge u'$$ in $V(s\oplus i)$, where
$u'=u$ or $u'=0$. If $x_i\ne 0\in V(s\oplus i)$ then $x_i=0\in V(s)$ and 
 $$\partial_{s\oplus i}^s(r^s_j\wedge
u)=x_i\wedge r^s_j\wedge u=x_i\wedge r^{s\oplus i}_j\wedge u+\alpha_j
x_i\wedge x_i\wedge u= r^{s\oplus i}_j \wedge (x_i\wedge u).$$ In
any case the map $\partial_{s\oplus i}^s$ is well defined.

\medskip

\noindent{\em Commutativity of 2-faces.}

Every 2-face of hypercube looks like
$$
\xymatrix{ \bigwedge^*V(s\oplus j)\ar[r]^{\partial^{s\oplus j}_{s\oplus i\oplus j}} & \bigwedge^*V(s\oplus i\oplus j)\\
\bigwedge^*V(s)\ar[r]_{\partial^{s}_{s\oplus
i}}\ar[u]^{\partial^{s}_{s\oplus j}} & \bigwedge^*V(s\oplus
i).\ar[u]_{\partial^{s\oplus i}_{s\oplus i\oplus j}} }
$$

According to dimensions of spaces $V(s'),\ s'=s, s\oplus i, s\oplus
j, s\oplus i\oplus j$ we have five types of diagrams without odd
edges:

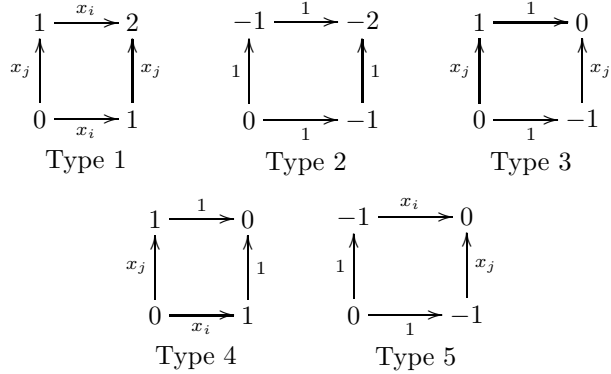
\begin{gather*}
\tar{\xymatrix{ 1\ar[r]^{x_i} & 2\\
0\ar[r]_{x_i}\ar[u]^{x_j} & 1\ar[u]_{x_j}}}{Type 1}\quad
\tar{\xymatrix{ -1\ar[r]^{1} & -2\\
0\ar[r]_{1}\ar[u]^{1} & -1\ar[u]_{1}}}{Type 2}\quad
\tar{\xymatrix{ 1\ar[r]^{1} & 0\\
0\ar[r]_{1}\ar[u]^{x_j} & -1\ar[u]_{x_j}}}{Type 3}\\
\tar{\xymatrix{ 1\ar[r]^{1} & 0\\
0\ar[r]_{x_i}\ar[u]^{x_j} & 1\ar[u]_{1}}}{Type 4}\quad
\tar{\xymatrix{ -1\ar[r]^{x_i} & 0\\
0\ar[r]_{1}\ar[u]^{1} & -1\ar[u]_{x_j}}}{Type 5}
\end{gather*}

Here the number at the place of the state $s'$ is equal to $\dim
V(s')-\dim V(s)=\cor A(s')-\cor A(s)$ and the label $z=1, x_i, x_j$
at the arrow for the map $\partial^{s'}_{s''}$ means that
$\partial^{s'}_{s''}(u) =z\wedge u$.

2-faces of types 1,2,3 are obviously commutative. Any 2-face of type
4 is commutative because $x_i=x_j=0\in V(s\oplus i\oplus j)$. For
the face of type 5 we need to show that $x_i=x_j\in V(s\oplus
i\oplus j)$. There are three possibilities.

1. $\sgn(v_i)=\sgn(v_j)=-1$. Then $v_i, v_j\in s\oplus i\oplus j$.
Without loss of generality we can assume that $v_i$ and $v_j$ are
the last vertices in $s\oplus i\oplus j$. The adjacency matrix
$A(s\oplus i\oplus j)$ can be represented in the form
\begin{equation}\label{Aab}
\bega{ccc}
A & a^\top & b^\top\\
a & \alpha & {\gamma}\\
b & {\gamma} & \beta\ena
\end{equation}
where $A=A(s)$.

Since $\cor A(s\oplus i) =\cor A(s)-1$ we have $x_i=0\in V(s\oplus
i)$. Then the rows of the matrix $A(s\oplus i)=\bega{cc} A& a^\top\\ a&
\alpha\ena$ generate the vector $(0\ 1)$. Hence, the rows of the
matrix $A(s\oplus i\oplus j)$ generate the vector $(0\ 1\ \delta),\
\delta\in\Z_2$. If $\delta=0$ then $x_i=0\in V(s\oplus i\oplus j)$
but this is not the case. So $\delta=1$ and we have the relation
$x_i+x_j=0$ in $V(s\oplus i\oplus j)$.

2. $\sgn(v_i)=-1,\ \sgn(v_j)=1$. Then $v_i\in s\oplus i\oplus j$ and
$v_j\not\in s\oplus i\oplus j$. Assume that the adjacency matrix
$A(s\oplus i)$ has the form~(\ref{Aab}) where $A=A(s\oplus j).$

We have the equality $x_i=0\in A(s\oplus i)$. It implies that the vector
$(0\ 1\ 0)$ is a linear combination of the rows of the adjacency
matrix $A(s\oplus i)$. If the coefficient in the linear combination
by the row $(b\ \gamma\ \beta)$ is zero then the rows of the matrix
$A(s\oplus i\oplus j) = \bega{cc} A& a^\top\\ a & \alpha\ena$
generate the vector $(0\ 1)$. Thus $x_i=0\in V(s\oplus i\oplus j)$
that leads to contradiction. Therefore the coefficient by $(b\
\gamma\ \beta)$ is $1$. Then the vector $(b\ \gamma+1)$ is generated
by the rows of the matrix $A(s\oplus i\oplus j)$ that means
$x_i+x_j=0\in V(s\oplus i\oplus j)$.

3. $\sgn(v_i)=\sgn(v_j)=-1$. Then the matrix $A(s)$ looks
like~(\ref{Aab}) where $A=A(s\oplus i\oplus j)$. The equality
$x_i+x_j=0\in V(s\oplus i\oplus j)$ corresponds to the relation
vector
$a+b$. Since $\rank A(s)=\rank A(s\oplus j) = \rank\bega{cc} A & a^\top\\
a & \alpha\ena = \rank\bega{ccc} A & a^\top& b^\top\\
a & \alpha & \gamma \ena$ the last row in the matrix $A(s)$ is
expressed as a linear combination of the other rows. If the
coefficient by the row $(a\ \alpha\ \gamma)$ in the combination is
zero then the vector $(b\ \gamma\ \beta)$ depends on the rows of the
matrix $\bega{ccc} A & a^\top & b^\top\ena$ and the vector $b$
depends on the rows of the matrix $A$ . Hence, $x_j=0\in V(s\oplus
i\oplus j)$ that is not true. Thus, the coefficients by $(a\ \alpha\
\gamma)$ is not zero so the vector $a+b$ is generated by the rows of
the matrix $A$ that means $x_i+x_j=0\in V(s\oplus i\oplus j)$.

Any diagram, which has odd edges in both upper-right and left-upper paths,
is commutative because the map of any odd edge is zero. After this
remark it remains only two nontrivial diagrams (up to symmetry between $i$ and
$j$) with an odd edge:

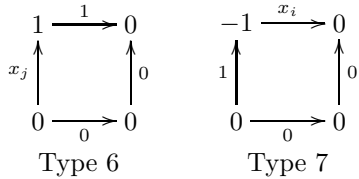
\begin{gather*}
\tar{\xymatrix{ 1\ar[r]^{1} & 0\\
0\ar[r]_{0}\ar[u]^{x_j} & 0\ar[u]_{0}}}{Type 6}\quad
\tar{\xymatrix{ -1\ar[r]^{x_i} & 0\\
0\ar[r]_{0}\ar[u]^{1} & 0\ar[u]_{0}}}{Type 7}
\end{gather*}

The diagram of type 6 is commutative because $\dim V(s\oplus i)=\dim
V(s\oplus i\oplus j)$ so $x_j=0\in V(s\oplus i\oplus j)$.

For the commutativity of the diagram of type 7 we need to prove that
$x_i= 0\in V(s\oplus i\oplus j)$. Assume this is not the case.

%1.
Let $\sgn(v_i)=\sgn(v_j)=-1$. Then the adjacency matrix looks
like~(\ref{Aab}). Since $x_j= 0\in V(s\oplus i)$ we have
$\rank\bega{cc} A &a^\top\\ a& \alpha \\ b & \gamma \ena =
\rank\bega{cc} A &a^\top\\ a& \alpha\ena$. Then $b$ depends on the
rows of $\bega{c} A \\ a\ena$. But $\rank \bega{c} A\\ a\ena = \rank
A$ hence $b$ depends on the rows of $A$. On the other hand, we have
$x_j\ne 0\in  V(s)$ so $\rank\bega{c}A\\b\ena\ne \rank A.$ This is
contradiction.

The other cases of sign arrangement can be treated in an analogous
manner.

%2. $\sgn(v_i)=\sgn(v_j)=-1$. Then we can assume $A(s\oplus i\oplus
%j)= \bega{cc} A &a^\top\\ a& \alpha$, $A(s\oplus j)= A$

This concludes the proof of Proposition~\ref{prop_val}.
\end{proof}

\begin{defi}
We call the homology $\Kh(G)$ of the complex $(C(G),\partial)$  the {\em
reduced (odd) Khovanov homology} of the labeled simple graph $G$.
\end{defi}

The main theorem of the article states that Khovanov homology is
well defined for graph-links.

\begin{theor}
Khovanov homology $\Kh(G)$ is invariant under $\Omega_1-\Omega'_4$
moves.
\end{theor}

\begin{proof}
Let $G$ be a labeled graph and $\wt G$ be the graph obtained from
$G$ by some Reidemeister move $\Omega_1-\Omega'_4$.

\medskip
\noindent{\em Invariance under $\Omega_1$.}

Let $\widetilde G$ be obtained from $G$ by an addition of an isolated
labeled vertex $v$. The complex $C(\wt G)$ is isomorphic to the
product of complexes $C(G)\otimes C(v)$ where the complex $C(v)$ is
equal to
$$\xymatrix{\Z_2 \ar[r]^-{x\wedge} & \bigwedge^*\Z_2\langle
x\rangle}$$ if $\sgn(v)=-1$ and
$$\xymatrix{\bigwedge^*\Z_2\langle x\rangle \ar[r]^-{x=0} & \Z_2}$$ if
$\sgn(v)=1$. In any case $H_*(C(v))=\Z_2\cdot 1$, where $1\in
H_0(C(v))$ if $\sgn(v)=1$ and $1\in H_1(C(v))$ if $\sgn(v)=-1$.
Thus, we have
$$ \Kh(\wt G) = \Kh(G)\otimes \Kh(v)\cong\Kh(G).$$

\noindent{\em Invariance under $\Omega_2$.}

Assume that we add the vertices $v$ and $w$ to get the graph $\wt G$
by $\Omega_2$ and $\sgn(v)=1$, $\sgn(w)=-1$ . Without loss of
generality we can write the adjacency matrix $A(\wt G)$ in one of
this two forms
$$
\bega{ccc}  0&0& a^\top \\
0&0& a^\top\\
a&a&A(G)\ena \quad\mbox{ or }\quad
\bega{ccc}  1&1& a^\top \\
1&1& a^\top\\
a&a&A(G)\ena.
$$
In both cases for every state $s\in \V(G)$ we have the following
equalities: $\cor A(G(s))=\cor A(\wt G(s))$, $\cor A(\wt
G(s\cup\{v\}))=\cor A(\wt G(s\cup\{w\}))=\cor A(\wt
G(s\cup\{v,w\}))-1$. These equalities define the type of the upper
and left arrows of the complex $C(\wt G)$ written in the form
$$
\xymatrix{C_{vw}\ar[r]^{1} & C_w\\
C_v\ar[u]^{x_2}\ar[r]_\partial & C.\ar[u]_\partial}
$$
Here $C_v$ consists of the chains whose state contains $v$ and
does not contain $w$; the complexes $C, C_w, C_{vw}$ are defined analogously.

For any state $s$ in $C_{vw}$ let us define a linear function $f :
V(s)\to \Z_2$ by the formula
$f(\sum_i\lambda_ix_i)=\lambda_1+\lambda_2$. The function $f$ is
well defined because it equals to zero on any relation:
$f(r^s_i)=a_{i1}+a_{i2}=0$ since $a_{i1}=a_{i2}$. Then
$\bigwedge^*V(s)=\bigwedge^*\ker f\oplus x_2\bigwedge^*V(s)$ and
$C_{vw}=X\oplus x_2C_{vw}$. The subcomplex
$X\rightarrow C_{w}$ is acyclic because the arrow is an isomorphism. Then the homology of $C(\wt G)$
coincides with the homology of the quotient complex
$$
\xymatrix{x_2C_{vw} & \\
C_v\ar[u]^{x_2}\ar[r]_\partial & C.}
$$
The quotient of the obtained complex by the subcomplex $C$ appears
to be acyclic too. Thus $C(\wt G)$ has the same homology as
$C=C(G)$.

\medskip
\noindent{\em Invariance under $\Omega_3$.}

Without loss of generality we can assume that the vertices $u, v, w$
in the third Reidemeister move have the indices $1, 2, 3$ in 
$\V(G)=\V(\wt G)$. There are two variants of $\Omega_3$ depending
adjacency of $v$ and $w$ and we consider the version where $v$ and
$w$ are not adjacent. The proof for the other version is analogous.
The adjacency matrices of $G$ and $\wt G$ looks like
$$
A(G)=\bega{cccc} 0&1&1&0^\top\\
1&0&0&a^\top\\
1&0&0&b^\top\\
0&a&b&B\ena,\quad
A(\wt G)=\bega{cccc} 0&0&0&(a+b)^\top\\
0&0&0&a^\top\\
0&0&0&b^\top\\
a+b&a&b&B\ena.
$$

Denote $\wt V(s)=V(\wt G(s))$. Then for any $s\subset
\V(G)\setminus\{u,v,w\}$ we have $V(s)\cong \wt V(s)$, $V(s\oplus
v)\cong \wt V(s\oplus v)$, $V(s\oplus w)\cong \wt V(s\oplus w)$,
$V(s\oplus v\oplus w)\cong \wt V(s\oplus u\oplus v)\cong \wt
V(s\oplus u\oplus w)$, $V(s\oplus u\oplus v\oplus w)\cong \wt
V(s\oplus u)$, $V(s\oplus u\oplus v)\cong V(s\oplus u\oplus w)\cong
V(s)$ and the correspondent isomorphisms of the exterior algebras are
compatible with the differential.

Consider complexes $C(G)$ and $C(\wt G)$ in the form of cubes:
$$
\xymatrix@!0{ & C_{uw} \ar[rr] & & C_{uvw} \\
C_u \ar[ur]^1\ar[rr]^(0.65)1 & & C_{uv}\ar[ur] &\\
& C_w\ar'[r][rr]\ar'[u][uu] && C_{vw}\ar[uu]\\
C\ar[rr]\ar[ur]\ar[uu]^{x_1} && C_v\ar[ur]\ar[uu] &} \qquad
\xymatrix@!0{ & \wt C_{uv} \ar[rr] & & \wt C_u \\
\wt C_{uvw} \ar[ur]^1\ar[rr]^(0.65)1 & & \wt C_{uw}\ar[ur] &\\
& \wt C_v\ar'[r][rr]\ar'[u][uu] && \wt C\ar[uu]\\
\wt C_{vw}\ar[rr]\ar[ur]\ar[uu]^{x_1} && \wt C_w\ar[ur]\ar[uu] &}
$$

For any state $s$ in $C_{u}$ let us define a linear function $f :
V(s)\to \Z_2$ by the formula $f(\sum_i\lambda_ix_i)=\lambda_1$. The
function is well defined and there are decompositions
$\bigwedge^*V(s)=\bigwedge^*\ker f\oplus x_1\bigwedge^*V(s)$ and
$C_{u}=X\oplus x_1C_{u}$. Consider the following subcomplex
$$
\xymatrix@!0{ & C_{uw} \ar[rr] & & C_{uvw} \\
X \ar[ur]^1\ar[rr]^(0.65)1 & & C_{uv}\ar[ur] &\\
& C_w\ar'[r][rr]\ar'[u][uu] && C_{vw}\ar[uu]\\
 && C_v\ar[ur]\ar[uu] &}
$$
The quotient complex $C\to x_1C_{u}$ is acyclic so the homology of
$C(G)$ is isomorphic to the homology of the subcomplex. This
subcomplex contains the acyclic subcomplex $X\to \partial(X)$. The
maps $X\to C_{uv}$ and $X\to C_{uw}$ are isomorphisms, so after
factorization we obtain the complex
$$
\xymatrix@!0{ & C_{uw} \ar[rr]\ar@{=}[dr] & & C_{uvw} \\
& & C_{uv}\ar[ur] &\\
& C_w\ar'[r][rr]\ar[uu] && C_{vw}\ar[uu]\\
 && C_v\ar[ur]\ar[uu] &}
$$
where the spaces $C_v$ and $C_w$ are identified.

The reasonings analogous to the reasonings above reduce $C(\wt G)$
to the complex (for a state $s$ in $\wt C_{uvw}$ we should define
the function $f : \wt V(s)\to \Z_2$ by the formula
$f(\sum_i\lambda_ix_i)=\lambda_1+\lambda_2+\lambda_3$).
$$
\xymatrix@!0{ & \wt C_{uv} \ar[rr]\ar@{=}[dr] & & \wt C_{u} \\
& & \wt C_{uw}\ar[ur] &\\
& \wt C_v\ar'[r][rr]\ar[uu] && C\ar[uu]\\
 && \wt C_w\ar[ur]\ar[uu] &}
$$
Both complexes are isomorphic to the complex $C_v\oplus C_w\to
C\oplus C_{vw}\to C_{uvw}$. Thus $C(G)$ and $C(\wt G)$ have the same
homology.

\medskip
\noindent{\em Invariance under $\Omega_4$.}

Let the vertices $u$ and $v$ of the move $\Omega_4$ have the indices $p$
and $q$ in $V(G)=V(\wt G)$. The coefficients of adjacency matrices
of $A(G)=(a_{ij})$ and $A(\wt G)=(\wt a_{ij})$ are connected by the
formula
$$
\wt a_{ij} = \left\{\begin{array}{cl}
a_{ij}+a_{ip}a_{jq}+a_{iq}a_{jp},& \{i,j\}\cap \{p,q\}=\emptyset,\\
a_{ij}, & \{i,j\}\cap \{p,q\}\ne\emptyset.
\end{array}\right.
$$

Consider the map $\phi$ acting on the states by
the formula
$$
\phi(s)=\left\{\begin{array}{cl} s\cup\{u,v\},& \{u,v\}\cap
s=\emptyset, \\
s\setminus \{u,v\}, & \{u,v\}\cap s=\{u,v\},\\
s, & \{u,v\}\cap s\ne\emptyset, \{u,v\}\\
\end{array}\right.
$$
together with the linear maps $\Phi : V(s)\to \wt V(\phi(s))$ defined by the
formula

$$
\Phi(x_i)=\left\{\begin{array}{cl} x_i,& i\ne p,q, \\
x_q, & i=p,\\
x_p, & i=q.\\
\end{array}\right.
$$

One can check that the map $\Phi$ defines an isomophism of the
vector spaces and after the natural extension to an isomorphism of
external algebras $\bigwedge^*V(s)\to\bigwedge^*\wt V(\phi(s))$ for each state $s$ it
determines an isomophism of the graded linear spaces $\Phi:C(G)\to
C(\wt G)$. $\Phi$ is a chain map. Thus the complexes $C(G)$ and
$C(\wt G)$ are isomorphic as well as their homology.

\medskip
\noindent{\em Invariance under $\Omega'_4$.}

Let the vertex $v$ of the move $\Omega'_4$ have the index $p$. The
coefficients of adjacency matrices of $A(G)=(a_{ij})$ and $A(\wt
G)=(\wt a_{ij})$ are connected by the formula
$$
\wt a_{ij} = \left\{\begin{array}{cl}
a_{ij}+a_{ip}a_{jp},& i,j\ne p,\\
a_{ip}, & j=p,\\
a_{pj}, & i=p,\\
\end{array}\right.
$$

Consider the map $\phi: C(G)\to C(\wt G)$ which acts on the states
by the formula
$$
\phi(s)=s\oplus p
$$
and the linear maps $\Phi : V(s)\to \wt V(\phi(s))$ defined on the generators by the
formula $\Phi(x_i)=x_i,\ 1\le i\le n$. Then $\Phi$ induces a well
defined isomorphism of the complexes $C(G)$ and $C(\wt G)$ so the
homologies of the complexes coincide.
\end{proof}

\begin{coro}
Khovanov homology $\Kh(G)$ is an invariant of graph-links.\hfill
$\Box$
\end{coro}

\section{Jones polynomial of graph-knots and Khovanov homology}

Let $G$ be a simple labeled graph. The complex $C(G)$ has two
grading: the homological grading $M_0$ and the algebraic grading
$\deg$ of graded algebras $\bigwedge^*V(s)$. The differential is
not homogeneous with respect to $\deg$ but it is compatible with the
grading $Q_0$ where for element any $u\in \bigwedge^r (s)$ we define
$Q_0(u)=\dim_{\Z_2}V(s)-2r+M_0(s)$. The differential increases the
grading $M_0$ and leaves the grading $Q_0$ unchanged. Then
$$
\Kh(G)=\bigoplus_{m,q\in\Z}\Kh(G)_{(m,q)}.
$$

The following definition is due to Ilyutko and Manturov~\cite{IM}.

\begin{defi}
Let $G$ be a simple labeled graph with $n$ vertices. The {\em
Kauffman bracket polynomial} of $G$ is defined as
$$
\langle G\rangle(a) = \sum_{s\subset
\V(G)}a^{\alpha(s)-\beta(s)}(-a^2-a^{-2})^{\cor A(s)},
$$
where $\alpha(s)=\#\{ v\in s\,|\,\sgn(v)=-1\}+\#\{ v\not\in
s\,|\,\sgn(v)=1\}$ and $\beta(s)=n-\alpha(s)$.
\end{defi}

It appears that the Poincar\'e polynomial of the bigraded Khovanov
homology of the graph $G$ coincides with its Kauffman bracket up to
shift of the gradings.
\begin{theor}\label{kauff}
$$ \sum_{m,q\in\Z}(-1)^m\dim_{\Z_2}\Kh(G)_{(m,q)}\cdot t^q=(it^{1/2})^{-n}\langle G\rangle(it^{1/2}).$$
\end{theor}
\begin{proof}
The left term of the equality coincides with the Poincar\'e
polynomial of the chain complex $C(G)$. Each state $s\subset\V(G)$
has homological grading $\alpha(s)$ and the corresponding chain
space $\bigwedge^*V(s)$ contributes
$(-1)^{\alpha(s)}(t+t^{-1})^{\dim V(s)}t^{\alpha(s)}$ to the
Poincar\'e polynomial. Since $\dim V(s)=\cor A(s)$ and
$\alpha(s)=\frac 12 (\alpha(s)-\beta(s)-n)$ the sum over all states
yields the right term of the equality.
\end{proof}

The gradings $M_0, Q_0$ of the Khovanov complex are not preserved by
the Reidemeister moves as shown in the following table.
\begin{equation}\label{grad}
\begin{array}{|c|c|c|}
\hline & M_0 & Q_0 \\
\hline \Omega_1^+ & 0 &-1\\
\Omega_1^- & 1 &1\\
\Omega_2 & 1 &0\\
\Omega_3 & 0 &0\\
\Omega_4 & 0 &0\\
\Omega'_4 & 0 &0\\
\hline
\end{array}
\end{equation}
Here $\Omega_1^\pm$ denotes the adding an isolated vertex with the
label $\pm 1$ and $\Omega_2$ denotes the adding two vertices. The
entries of the cells are the shifts of the gradings after the
corresponding move.

 Thus the groups $\Kh(G)_{(m,q)}$ are not invariants of
graph-links. Nevertheless we can normalize the gradings for the
graph-knots.

\begin{defi}[\cite{IM}] A graph-link $G$ is called a {\em graph-knots}
if $\cor(A(G)+E)=0$ where $E$ is the identity matrix.
\end{defi}
The condition  $\cor(A(G)+E)=0$ survives after Reidemeister moves so
the definition is correct.

\begin{defi}[\cite{IM}]
{\em Writhe number} of a graph $G$ with $n$ vertices is the sum
$$
w(G)=\sum_{i=1}^n(-1)^{\cor B_i(G)}\sgn v_i,
$$
where $B_i(G)=A(G)+E+E_{ii}$ and the only non-zero element of
$E_{ii}$  lies in the $i$-th row and $i$-th column.
\end{defi}

Writhe number is invariant under the moves $\Omega_2-\Omega_4$. The
move $\Omega_1^\pm$ changes the writhe number by $\mp 1$. Using this
fact and the table~(\ref{grad}) we construct two gradings which are
invariant under the Reidemeister moves
$$ M = M_0-\frac{n(G)+w(G)}2,\quad Q=Q_0+\frac{n(G)+3w(G)}2
$$
where $n(G)$ is the number of vertices in $G$.

Let us denote $\Kh_{m,q}(G)$ to be the homogeneous part of $\Kh(G)$ with
the gradings $M=m$ and $Q=q$. The groups $\Kh_{m,q}(G)$ are
invariants of graph-knots.

Writhe number allows to introduce another invariant of graph-knots.
\begin{defi}[\cite{IM}]
{\em Jones polynomial} of a graph-knot $G$ is defined as
$$
X(G)(a)=(-a)^{-3w(G)}\langle G\rangle(a).
$$
\end{defi}

Theorem~\ref{kauff} leads to the following statement.

\begin{theor}
The bigraded Khovanov homology of graph-knots $\Kh_{m,q}(G)$
categorifies the Jones polynomial in the sense that
$$ \sum_{m,q\in\Z}(-1)^m\dim_{\Z_2}\Kh_{m,q}(G)\cdot t^q= X(G)(it^{1/2}).$$
\end{theor}

The author is grateful to V.~Manturov and D.~Ilyutko for useful
discussions.

\end{document}